\input amstex

\def\a{\alpha}
\def\b{\beta}
\def\l{\lambda}

\def\ee#1{e_{#1}}

\def\hd{, \hdots ,}

\def\na{n+a}
\def\ooo#1#2{\omega^{#1}_{#2}}
\def\oo#1{\omega^{#1}}

\def\pp#1{\Bbb P^{#1}}
\def\ppp{\Bbb P}

\def\tdim{\text{dim}\,}

\def\tmod{\text{ mod }}

\documentstyle{amsppt}
\magnification = 1100
\hsize =15truecm
\hcorrection{.5truein}
\baselineskip =18truept
\vsize =22truecm
\NoBlackBoxes
\topmatter
\title Lines on projective varieties
\endtitle
  \author
    J.M. Landsberg 
  \endauthor


\address{School of Mathematics,
Georgia Institue of Technology,
 686 Cherry St., Skiles Bldg.
 Atlanta, GA 30332-0160
 USA  }\endaddress
\email {jml\@math.gatech.edu }
\endemail
\endkeywords
\subjclass{ primary 53, secondary 14}\endsubjclass
 
\endtopmatter

\document

\subheading{Definition} 
A variety (or manifold) $X\subset\pp N$ (or $\Bbb A^N$)
is {\it covered by lines} if through a general point $x\in X$ there
passes a finite number of lines contained in $X$. 

It was known classically
that a surface covered by lines can contain at most two lines through
a general point.
 In [MP] it was shown that a
$3$-fold covered by lines contains at most $6$ lines through
a general point. 

I prove:

\proclaim{Theorem 1} Let $X^n\subset \pp\na$ be covered 
by lines. Then there are at most $n!$ lines passing through a general
point of $X$.
\endproclaim

The theorem holds for any ground field
of characteristic zero. The same conclusion is valid for analytic  varieties
in affine or projective space if one asks that the lines be
  contained in the closure of $X$. The same conclusion even holds
in the $C^{\infty}$ category if one replaces \lq\lq general point\rq\rq\
by \lq\lq every point\rq\rq\ in the hypotheses.

Theorem 1 will be a consequence of theorem 2:

\proclaim{Theorem 2} Let $X^n\subset\pp{n+1}$ be a hypersurface and
let $x\in X$ be a general point. Let $\Sigma^{\l}\subset \ppp T_xX$
denote the tangent directions to lines having contact to
order $\l$ with $X$ at $x$. (The notation is such that
$\Sigma^1=\ppp T_xX$.) If there is an irreducible component
$\Sigma^k_0\subset \Sigma^k$ such that $\tdim \Sigma^k_0 >n-k$ then
$\Sigma^k_0\subset\Sigma^{\infty}$, i.e., all lines corresponding to
points of
$\Sigma^k_0$ are contained in $X$.
\endproclaim

Note that the expected dimensionof $\Sigma^k$ is $n-k$.

\proclaim{Corollary}Let $X^n\subset \pp\na$ be a  variety
such that through a general point $x\in X$ there is a $p$-plane
 having contact to order $n-p+2$. Then $X$ is uniruled by $\pp p$'s.
\endproclaim

The corollary is not expected to be optimal for most
values of $k$, see [L2, theorem 4].

\demo{Proof of theorem 2}
We use notation as in [L1] and [L2]. 
 We choose a
basis $\ee 1\hd \ee n$ of $T_xX$
 such   $[\ee 1]$ is
a general point of  $\Sigma^k_0$,
and   $\tilde T_{[\ee 1]}\Sigma^k=\ppp \{ \ee 1,\ee 2\hd\ee p\}$.
By hypothesis $p-1=\tdim \Sigma^k_0 >n-k$.

Let $1\leq \a,\b\leq n$.
We let $r_{\a ,\b}$ denote the coefficients of the second fundamental form $F_2$
of $X$ at $x$ and $r_{\a_1\hd \a_i}$ denote the coefficients of $F_i$.
sometimes we write $r^i_{\a_1\hd \a_i}=r_{\a_1\hd \a_i}$ to 
keep track of   $i$.

We let $2\leq s,t\leq p$, $p+1\leq  j,l \leq n$.

By our normalizations,
we have $r^\l_{1\hd 1}=0$ and $r^\l_{1\hd 1,s}=0$ for
$2\leq s\leq p$ and  $2\leq \l\leq k$.
We will show that $r^h_{1\hd 1}=0$ and $r^h_{1\hd 1,s}=0$ for $2\leq s\leq p$
and 
for all $h$, showing that there is a $p$-dimensional space of lines passing
through $X$ at $x$.

The technique of proof is the same as in [L2], namely we will use
our lower order equations to solve for the connection forms
$\ooo j 1$ in terms of the semi-basic forms $\oo l$ and
then plug into the higher  $F_h$ to obtain the vanishing of  
$F_h(\ee 1\hd \ee 1)$ and $F_h(\ee 1\hd \ee 1,\ee s)$.
Using the formalism for
the $F_h$ developed in  from [L1],  we obtain that for $\l\leq k $:
 
$$
r^\l_{1\hd 1,i}\oo i  = -r^{\l-1}_{1\hd 1,j}\ooo j 1
$$
These are $k-2$ equations for the
$n-p$ one-forms $\ooo j 1$ in terms of the $n-p$ semi-basic forms
$\oo j$. Recall that $n-p\leq k-2$.
 The system is solvable  as
were the $(k-2)\times(n-p)$ matrix
$(r^\l_{1\hd 1,j})$ not of maximal rank, there would be additional directions
in the tangent space to $\Sigma^k_0$. Thus we have 
$$
\ooo j 1\equiv 0\tmod \{\oo {p+1}\hd\oo n\}
$$
so the equation
$$
r^{k+1}_{1\hd 1,\beta}\oo\beta  = -r^{k}_{1\hd 1,j}\ooo j1
$$
implies $r^{k+1}_{1\hd 1}=r^{k+1}_{1\hd 1,2}=\hdots
r^{k+1}_{1\hd 1,p}=0$ and thus the
line through $[\ee 1]$ has contact to order $k+1$
and moreover   $\tilde T_{[\ee 1]}
\Sigma^{k+1}=\tilde T_{[\ee 1]}\Sigma^{k }$.
  Now one can use
these equations iteratively to show the same holds to order $k+2$
and all orders, i.e., the component of $\Sigma^k$ containing $[\ee 1]$
equals the component of $\Sigma^{\infty}$ containing $[\ee 1]$.
\qed\enddemo

\demo{Proof of theorem 1} Without loss of generality in
theorem 1 we may assume $X$ is a hypersurface  as one
can reduce to this case by linear projection.
First note that $n!$ is the expected bound in the sense that in
$\ppp T_xX$ one has the ideal $I$ generated by $F_2,F_3\hdots F_{n }$
defining the variety $\Sigma^{n }\subset\ppp T_xX$ of all lines
having contact with $X$ at $x$ to order $n $.
(Note that the polynomials of degree greater than
two are not well defined individually
but the ideal $I$ is.) Since
$\ppp T_xX$ is a $\pp{n-1}$, if $\Sigma^{n }$
is zero-dimensional, it is at most $n!$ points and we  are done.
 If  $\tdim\Sigma^{n }>0$ then theorem 2 applies.\qed\enddemo

 \subheading{Concluding remarks}
 Theorem 1
 generalizes theorem 1 of [L2], which
 says that if $\Sigma^{n+1}\neq\emptyset$, then $\Sigma^{n+1}=
\Sigma^{\infty}$. It is sharp for hypersurfaces as a general hypersurface
 of degree $n$ will be covered by $n!$ lines.
   It is unlikely theorem 1 is
 sharp in higher codimension, however   the following example
  due to F. Zak shows that one cannot hope
  to do too much better: Let $Y^{n-1}$ be a hypersurface with $(n-1)!$ lines
 passing through a general point and let $X=Seg(Y\times C)$ where
 $C\subset\pp M$ is a curve. One can arrange for the codimension of $X$
 to be arbitrarily large in this way and $X$ has $(n-1)!$ lines passing
 through a general point. 
 
 One could hope to now do a classification in the spirit of [MP].
 For example, $n!$ should only be possible for hypersurfaces of degree
 $n$. An interesting class of examples with a small number of lines
 is obtained by taking linear sections of uniruled homogeneous varieties.

 \subheading{Acknowledgements} It is a pleasure to thank E. Mezzetti and
 F. Zak for useful conversations.

\enddocument